\documentclass[11pt]{amsart}
\usepackage{mathrsfs}
\usepackage{amssymb}
\usepackage{amsthm}
\usepackage{titletoc}
\usepackage{amssymb}
\usepackage{amsxtra}
\usepackage{mathrsfs}
\usepackage{amsmath,amstext,amsfonts,verbatim, hyperref}

\newtheorem{theorem}{Theorem}
\newtheorem{corollary}{Corollary}

\newtheorem{remark}{Remark}
\newtheorem{example}{Example}
\newtheorem{thmx}{Theorem}

\title{Difference independence of the Euler gamma function}
\author{Qiongyan Wang$^1$, Xiao Yao$^2$*}
\address{$^1$ School of Mathematical sciences, Peking University, Beijing, 100871, P.R. China
\vskip 2pt $^2$ School of Mathematical sciences and LPMC, Nankai University, Tianjin, 300071, P.R. China
\vskip 2pt \hspace{1.5mm} Email:{\sf qiongyanwang@aliyun.com, yaoxiao@nankai.edu.cn}}
\thanks{{\sf Corresponding author: Xiao Yao}}
\thanks{{\sf 2010 Mathematics Subject Classification.} 11M06, 39A05.}
\thanks{The research was partially supported  by National Key R\&D Program of China (2020YFA0713300) and NSFC of China(No.11901311).}
\thanks{{\sf Keywords.} Algebraic difference independence;  Euler gamma function; Algebraic difference equations.}

\begin{document}
\maketitle
\begin{abstract}
In this paper, we established a  sharp version  of the difference analogue of the celebrated H\"{o}lder's theorem concerning the differential independence of the Euler gamma function $\Gamma$.  More precisely, if $P$ is a polynomial of $n+1$ variables in $\mathbb{C}[X, Y_0,\dots, Y_{n-1}]$ such that
\begin{equation*}
P(s, \Gamma(s+a_0), \dots, \Gamma(s+a_{n-1}))\equiv 0
\end{equation*}
 for some $(a_0, \dots, a_{n-1})\in \mathbb{C}^{n}$ and $a_i-a_j\notin \mathbb{Z}$ for  any  $0\leq i<j\leq  n-1$, then we have $$P\equiv 0.$$

 Our result complements a  classical result of algebraic differential independence of the Euler gamma function  proved by H\"{o}lder  in 1886, and also a result of algebraic difference independence of the Riemann zeta function proved by Chiang and Feng in 2006.
\end{abstract}

\section{Introduction}
A classical theorem of H\"{o}lder \cite{6} states that the Euler gamma function
$$\Gamma(s)=\int_{0}^{+\infty}t^{s-1}e^{-t}dt, ~~~ \Re s> 0$$
which can be analytically continued to the whole complex plane $\mathbb{C}$,  does not satisfy any nontrivial algebraic differential equation whose coefficients are polynomials in $\mathbb{C}$. We state it in the following.
\begin{thmx}
Let $P$ be a polynomial of $n+1$ variables in $\mathbb{C}[X, Y_0,\dots, Y_{n-1}]$. Assume  that
\begin{equation*}
P(s, \Gamma(s), \dots, \Gamma^{(n-1)}(s))\equiv 0,
\end{equation*}
then we have
\begin{equation*}
P\equiv 0.
\end{equation*}
\end{thmx}

 To the best of our knowledge, the Euler gamma function $\Gamma$  seems to be the first known example which  satisfies the algebraic differential independence property in the literature. It is well known that the Riemann zeta function $\zeta$ is associated with $\Gamma$ by the famous Riemann functional equation
\begin{equation}\label{eq0.1}
\begin{aligned}
\zeta(1-s)=2^{1-s}\pi^{-s}\cos \frac{\pi s}{2}\Gamma(s)\zeta(s).
\end{aligned}
\end{equation}

Motivated by the Riemann functional equation, it is natural to consider the algebraic differential independence property for the Riemann zeta function. The study of the  algebraic differential independence of the Riemann zeta function $\zeta$ can be dated back to Hilbert.  In \cite{Hilbert}, he conjectured that H\"{o}lder's result can be extended to the Riemann zeta function $\zeta$. Later, this conjecture was verified by Ostrowski in \cite{11}.

 Bank and Kaufman  \cite{1, 2} made the following celebrated generalizations of  H\"{o}lder's result.
 \begin{thmx} \label{thm-Bank-Kaufman}
 Let $P$ be a polynomial in $K[X, Y_0, \dots, Y_{n-1}]$, where $K$ is the field of all meromorphic  functions   such that the Nevanlinna's characteristic
$T(r, f)=o(r)$ as $r$ goes to infinity for any $f$ in $K$. Assume that
\begin{equation*}
P(s, \Gamma(s), \dots, \Gamma^{(n-1)}(s))\equiv 0,
\end{equation*}
then we have
\begin{equation*}
P\equiv 0.
\end{equation*}
\end{thmx}

For the Nevanlinna characteristic $T(r,  f)$,  we refer to Hayman's book \cite{Hayman} for a detailed introduction. Since $\Gamma$ and $\zeta$ appeared very naturally in Riemann functional equation \eqref{eq0.1},  Markus in \cite{Markus} posted an open problem to study the joint algebraic differential independence of $\Gamma$ and $\zeta$.  We refer the readers to the references   \cite{7, 8, 8B, 9} for the recent  developments in this direction.

It is interesting to study the algebraic difference independence of $\zeta$ or $\Gamma$.  Feng and Chiang proved the following result.
\begin{thmx}\label{thm-Chiang-Feng}
Let $P$ be a polynomial of $n+1$ variables  in $\mathbb{C}[X, Y_0,\dots, Y_{n-1}]$ and $s_0, \dots, s_{n-1}$ be $n$ distinct  numbers in $\mathbb{C}$. Assume that
\begin{equation*}
P(s, \zeta(s+s_0), \dots, \zeta(s+s_{n-1}))\equiv 0,
\end{equation*}
then we have
\begin{equation*}
P\equiv 0.
\end{equation*}
\end{thmx}

 Chiang and Feng's result extended a result of Ostrowski  in \cite{11} where the assumption of $s_0, \dots s_{n-1}$ are $n$ distinct real numbers are needed.  Indeed, Chiang and Feng proved that Theorem \ref{thm-Chiang-Feng} also holds under  the same assumption in Theorem \ref{thm-Bank-Kaufman} , we refer the interested readers to \cite{3} for the details. Here, we also mention two remarkable universality results due to  Voronin in 1970s for the differential case \cite{Voronin-1975} and the difference case \cite{Voronin}. We refer to \cite{Steuding} for the detailed introduction of the recent developments in this direction.

To the best of our knowledge, the topic of the algebraic difference independence of the Euler gamma function was first  addressed by Hardouin in \cite{Hardouin} in the framework of difference Galois theory. Motivated by the multiplication theorem of Euler gamma function
\begin{equation}
\Gamma(ns)=n^{ns-\frac{1}{2}}(2\pi)^{\frac{1-n}{2}}\prod_{j=0}^{n-1}
\Gamma(s+\frac{j}{n}),
\end{equation}
Hardouin proved the following result.
\begin{thmx}[\cite{Hardouin}]\label{thm-D}
Let $a_0, \dots,  a_{n-1}$ be $n$ complex numbers in $\mathbb{C}$, $b_0, \dots, b_{m-1}(\geq 2)$ be $m$ integers such  that $\{a_j ( mod \ 1)\}_{j=0}^{n-1}$ and $\{\sum\limits_{l=0}^{b_j-1}\frac{l}{b_j}(mod \ 1) \}_{j=0}^{m-1}$ are $\mathbb{Z}$-linearly independent.  Assume that
\begin{equation*}
P(s, \Gamma(s+a_0), \dots, \Gamma(s+a_{n-1}), \Gamma(b_0s), \dots, \Gamma(b_{m-1} s))\equiv 0
\end{equation*}
for  some polynomial $P$, then we have
\begin{equation*}
P\equiv 0.
\end{equation*}
\end{thmx}
Hardouin's proof  relies on Kolchin's type theorem in an essential way. See also in \cite{BBD} for a detailed disccusion of Kolchin's type theorem and several powerful applications in  algebraic independence problems.

Our starting point is another well known difference equation of $\Gamma$,
\begin{equation}\label{eq-Gamma}
\Gamma(s+1)=s\Gamma(s).
\end{equation}
 This may be the obvious obstruction for us to study the  algebraic difference independence of the Euler gamma function $\Gamma$. One can not expect to obtain Theorem B  for $\Gamma$ directly.
While in this paper, we will show that  the machinery exhibited in  \eqref{eq-Gamma} is  the only obstruction to get the algebraic difference independence of $\Gamma$. Now, we state our main result in the following. In this paper, we will use an elementary method inspired by \cite{6, 11} to prove our main result, which avoids the advanced difference Galois theory. This may be of independent interest.

We define
\begin{equation}
\mathcal{H}:=\{(a_0, \dots, a_{n-1})\in \mathbb{C}^{n}: a_i-a_j\notin \mathbb{Z} \ \text{for any} \ 0\leq  i< j\leq n-1 \}.
\end{equation}
 Now, we state our main result in the following.
%

\begin{theorem}\label{thm-difference-Holder}
Let $P$ be a polynomial of $n+1$ variables in $\mathbb{C}[X, Y_0,\dots, Y_{n-1}]$.
 Assume that
\begin{equation*}
P(s, \Gamma(s+a_0), \dots, \Gamma(s+a_{n-1}))\equiv 0
\end{equation*}
for some $(a_0, \dots, a_{n-1})\in \mathcal{H}$,
then we have
\begin{equation*}
P\equiv 0.
\end{equation*}
\end{theorem}

We remark that we can also use Theorem \ref{thm-D} to recover part of the result of Theorem \ref{thm-difference-Holder} under the same condition of $(a _j)_{j=0}^{n-1}$ and also $m=0$ in Theorem \ref{thm-D}. While, it can not completely recover Theorem \ref{thm-difference-Holder}, since the condition  in Theorem \ref{thm-difference-Holder} is sharp.   Our result complements the classical result of algebraic differential independence of Euler gamma function  proved by H\"{o}lder \cite{6} in 1886, and also a result of algebraic difference independence of Riemann zeta function proved by Chiang and Feng \cite{3} in 2006.

\begin{corollary}
Let $P$ be a polynomial of $n+1$ variables in $\mathbb{C}[X, Y_0,\dots, Y_{n-1}]$. Assume that
\begin{equation*}
P(s, \Gamma(s), \dots, \Gamma(s+(n-1)\alpha))\equiv 0
\end{equation*}
for some $\alpha\not\in\mathbb{Q}$, then we have
\begin{equation*}
P\equiv 0.
\end{equation*}
\end{corollary}

\begin{remark}
 Theorem \ref{thm-difference-Holder} can be seen as a difference version of the H\"{o}lder's theorem. The identity \eqref{eq-Gamma} shows that the discussion restricted to $\mathcal{H}$   is  necessary.

We can also extend Theorem \ref{thm-difference-Holder} to the setting of $K[X, Y_0, \dots, Y_{n-1}]$ where $K$ is the field of all meromorphic  functions   such that the Nevanlinna's characteristic
$T(r, f)=o(r)$ as $r$ goes to infinity for any $f$ in $K$. While, we will not address it in this paper.
\end{remark}
 By Theorem \ref{thm-difference-Holder} and the Euclidean's algorithm, it is not hard to give the following two examples.
\begin{example}
Let $P=P(X, Y, Z)$ be a polynomial of $3$ variables in $\mathbb{C}[X, Y, Z]$. Assume that
\begin{equation*}
P(s, \Gamma(s+a_0), \Gamma(s+a_1))\equiv 0,
\end{equation*}
then
\begin{equation*}
P\equiv 0,
\end{equation*}
unless $a_1-a_0\in \mathbb{Z}$. In the latter case, if $\Re a_0<\Re a_1$,  $P$ can be divided by the polynomial $R(X, Y, Z)=Z-(X+a_0)\dots (X+a_1-1)Y$.
\end{example}

\begin{example}
Let $P(X, Y, Z, W)=YW-Z^2-YZ$ in $\mathbb{C}[X, Y, Z, W]$. We have
\begin{equation*}
P(s, \Gamma(s), \Gamma(s+1), \Gamma(s+2))\equiv 0.
\end{equation*}
$P$ belongs to the ideal
\begin{equation*}
<W-(X+1)Z, Z-XY>
\end{equation*}
generated by $W-(X+1)Z$ and $Z-XY$ in $\mathbb{C}[X, Y, Z, W]$. Furthermore, $P$ can be written by
\begin{equation*}
P(X, Y, Z, W)=Y(W-(X+1)Z)+Z(XY-Z).
\end{equation*}
\end{example}
\begin{remark}
Indeed, inspired by Example 1 and Example 2, we can apply Theorem \ref{thm-difference-Holder} and the Euclidean's algorithm again to give a complete characterization of the following set
\begin{equation*}
\mathcal{I}:=\{P\in \mathbb{C}[X, Y_0,\dots, Y_{n-1}]: P(s,\Gamma(s+a_0), \dots, \Gamma(s+a_{n-1}))\equiv 0\}
\end{equation*}
without  any assumption on $a_0, \dots, a_{n-1}$. While, we will not discuss it in this paper.
\end{remark}
\section{Proof of  Theorem \ref{thm-difference-Holder}}
 In order to prove Theorem \ref{thm-difference-Holder}, we need introduce a lexicographic order between any two monomials $Y_{0}^{i_0}\dots Y_{n-1}^{i_{n-1}}$ and $Y_{0}^{j_0}\dots Y_{n-1}^{j_{n-1}}$ in $\mathbb{C}[Y_0, \dots, Y_{n-1}]$, which plays an important role in our proof. And this strategty was inspired by Ostrowski's  proof of H\"{o}lder's classical proof  in \cite{11}. It also shares some  spirit of Kolchin's type theorem which was used in \cite{Hardouin, BBD}
 
We first introduce an order for the $n$ symbols $Y_0, \dots, Y_{n-1}$,
\begin{equation}\label{eq-order}
Y_0\prec Y_1\prec\dots\prec Y_{n-1}.
\end{equation}
This  can be used to induce a lexicographic order between any two monomials $Y_{0}^{i_0}\dots Y_{n-1}^{i_{n-1}}$ and $Y_{0}^{j_0}\dots Y_{n-1}^{j_{n-1}}$.  We still denote it by $\prec$ to simplify the notation. We define it in the following,
\begin{enumerate}
\item [{\bf case 1}:] $Y_{0}^{i_0}\dots Y_{n-1}^{i_{n-1}}=Y_{0}^{j_0}\dots Y_{n-1}^{j_{n-1}}$ if $i_k=j_k$ for $k=0, \dots, n-1$;
\item [{\bf case 2}:] $Y_{0}^{i_0}\dots Y_{n-1}^{i_{n-1}}\prec Y_{0}^{j_0}\dots Y_{n-1}^{j_{n-1}} $ if $i_0<j_0$ or
 there  exists  $1\leq k \leq n-1$ such that
 $$i_0=j_0, \dots, i_{k-1}=j_{k-1},  i_{k}<j_{k};$$
\item [{\bf case 3}:] $Y_{0}^{j_0}\dots Y_{n-1}^{j_{n-1}}\prec Y_{0}^{i_0}\dots Y_{n-1}^{i_{n-1}} $ can be defined similarly as in {\bf case 2}.
\end{enumerate}

For any nonzero polynomial $P=P(X, Y_0, \dots, Y_{n-1})$ in $\mathbb{C}[X, Y_0, \dots, Y_{n-1}]$, we  write it by
\begin{equation}\label{eq-sum-P}
P=\sum_{i=(i_0, \dots, i_{n-1})}\Phi_{i}(X)Y_{0}^{i_0}\dots Y_{n-1}^{i_{n-1}},
\end{equation}
where $\Phi_i(X)\in \mathbb{C}[X]$ and $\Phi_i(X)\neq  0$.
The {\bf highest term} of $P$ is defined by the  maximal element in
$\mathcal{T}_{P} $ with respect to the lexicographic order $\prec$ introduced above, where
\begin{equation}
\mathcal{T}_{P}:=\{Y_{0}^{i_0}\dots Y_{n-1}^{i_{n-1}}: \Phi_i(X)\ Y_{0}^{i_0}\dots Y_{n-1}^{i_{n-1}} \ \text{appeared in} \ \eqref{eq-sum-P}\}.
\end{equation}
 For any monomial $L=Y_{0}^{i_0}Y_{1}^{i_1}\dots Y_{n-1}^{i_{n-1}}$, we define its {\bf degree} $\deg(L)$  by $$\deg(L):=\sum_{k=0}^{n-1}i_{k}.$$
The {\bf height}  of $P$  is defined by the  degree of the highest term of $P$.

Now, we will prove Theorem \ref{thm-difference-Holder}.

\begin{proof}
Let
\begin{equation}
\mathcal{S}:=\{P\in \mathbb{C}[X, Y_0, \dots, Y_{n-1}]: P(s, \Gamma(s+a_0), \dots, \Gamma(s+a_{n-1}))\equiv 0\}.
\end{equation}
We will prove Theorem \ref{thm-difference-Holder} by contradiction. We assume that
$\mathcal{S} \neq \{0\}$. By our assumption, there exists a nonzero polynomial
\begin{equation*}
Q=\sum\limits_{i=(i_0, \dots, i_{n-1})}\Psi_{i}(X)Y_{0}^{i_0}\dots Y_{n-1}^{i_{n-1}},
\end{equation*}
which is of the lowest height in $\mathcal{S}\backslash\{0\}$ with $\Psi_j(X)Y_{0}^{j_0}\dots Y_{n-1}^{j_{n-1}}$ being its highest term for some $j=(j_0, \dots, j_{n-1})$.  Moreover, we also make the following assumption.

{\noindent \bf Assumption LD:} The  nonzero polynomial $\Psi_j(X)$ appearing in the highest term of $Q$ is also of the lowest degree.

Let
\begin{equation}\label{eq-def-T}
T(X, Y_0, \dots, Y_{n-1}):=Q(X+1, (X+a_0)Y_0, \dots, (X+a_{n-1})Y_{n-1}).
\end{equation}
Noting that
\begin{equation*}
Q(s, \Gamma(s+a_0), \dots, \Gamma(s+a_{n-1}))\equiv 0,
\end{equation*}
we have
$$
T(s, \Gamma(s+a_0), \dots, \Gamma(s+a_{n-1}))\equiv 0
$$
by \eqref{eq-Gamma}.
And the highest term of $T$ is $\hat{\Psi}_{j}(X)Y_{0}^{j_0}\dots Y_{n-1}^{j_{n-1}}$,
where
\begin{equation*}
\hat{\Psi}_{j}(X):=\Psi_{j}(X+1)(X
+a_0)^{j_0}\dots (X+a_{n-1})^{j_{n-1}}.
\end{equation*}

It follows from  the Euclidean's algorithm, there exist two polynomials $R=R(X)$ and $U=U(X)$ in $\mathbb{C}[X]$
such that
\begin{equation*}
\hat{\Psi}_j=R\Psi_j+U,
\end{equation*}
where either $U=0$ or $0<\deg U<\deg \Psi_j$. It is easy to see that $\deg R\geq 1$.

We claim that $U=0$. Otherwise, we know that the  polynomial
\begin{equation*}
H(X, Y_0, \dots, Y_{n-1}):=T(X, Y_0, \dots, Y_{n-1})-R(X)Q(X, Y_0, \dots, Y_{n-1})
\end{equation*}
is in $\mathcal{S}$.  It follows that the highest term of $H$ is
\begin{equation*}
U(X)Y_{0}^{j_0}\dots Y_{n-1}^{j_{n-1}}
\end{equation*}
 and $0<\deg U<\deg \Psi_j$. Thus, $H\neq 0$, which contradicts  the choice of $Q$ and  {\bf Assumption LD}. Now, we have
 $U=0$.

Since $U=0$, we  see that the highest term of $H$ is less than
the highest term of $Q$ if $H\neq 0$. This again contradicts our choice of $Q$. Thus, we get $H= 0$.  That is,
\begin{equation}\label{eq-TQ}
T(X, Y_0, \dots ,Y_{n-1})=R(X)Q(X, Y_0, \dots, Y_{n-1}).
\end{equation}

We first assume that  there exists $\beta\notin \Lambda :=\{-a_k:  0\leq k\leq n-1\}$ such that $R(\beta)=0$. By \eqref{eq-def-T} and \eqref{eq-TQ}, we  get
\begin{equation*}
Q(\beta+1, (\beta+a_0)Y_0, \dots, (\beta+a_{n-1})Y_{n-1})=0
\end{equation*}
in $\mathbb{C}[Y_0, \dots, Y_{n-1}]$. This implies that
\begin{equation*}
Q(\beta+1, Y_0, \dots, Y_{n-1})=\sum\limits_{i=(i_0, \dots, i_{n-1})}\Psi_{i}(\beta+1)Y_{0}^{i_0}\dots Y_{n-1}^{i_{n-1}}= 0
\end{equation*}
in $\mathbb{C}[Y_0, \dots, Y_{n-1}]$. Thus, we have
\begin{equation*}
\Psi_{i}(\beta+1)=0
\end{equation*}
for all $i$, which implies  that each $\Psi_i(X)$ can be divided by $X-\beta-1$. This contradicts  our assumption that $\Psi_j$ is of the lowest degree.

 Hence, each root of $R$ lies in $\Lambda$. Without loss of generality, we assume that $R(-a_0)=0$. Thus, we get
\begin{equation*}\label{eq-first-level}
Q(-a_0+1, 0, (a_1-a_0)Y_1, \dots, (a_{n-1}-a_0)Y_{n-1})= 0
\end{equation*}
by \eqref{eq-def-T} and \eqref{eq-TQ}.
Recalling that $a_j-a_0\notin \mathbb{Z}$ for any $j\neq 0$, we have
\begin{equation}\label{eq-Q}
Q(-a_0+1, 0, Y_1, \dots, Y_{n-1})= 0.
\end{equation}
Taking $X=-a_0+1$, $Y_0=0$ in \eqref{eq-def-T} and \eqref{eq-TQ}, we  get
\begin{align*}
&Q(-a_0+2, 0, (a_1-a_0+1)Y_1, \dots, (a_{n-1}-a_0+1)Y_{n-1})\\
=&R(-a_0+1)Q(-a_0+1, 0, Y_1, \dots, Y_{n-1})= 0
\end{align*}
by \eqref{eq-Q}.
Noting that  $a_j-a_0\notin \mathbb{Z}$ for any $j\neq 0$ again, we obtain
\begin{equation*}
Q(-a_0+2, 0, Y_1, \dots, Y_{n-1})=0
\end{equation*}
in $\mathbb{C}[Y_0, \dots, Y_{n-1}]$.
By induction, we can prove that for  any $m\in \mathbb{N}$,
\begin{equation*}
Q(-a_0+m, 0, Y_1, \dots, Y_{n-1})= 0
\end{equation*}
in $\mathbb{C}[Y_0, \dots, Y_{n-1}]$.
It follows by  the fundamental theorem of algebra, we get
\begin{equation*}
Q(X, 0, Y_1, \dots, Y_{n-1})=0
\end{equation*}
in $\mathbb{C}[X, Y_0, \dots, Y_{n-1}]$.
Thus, we proved that $Q$ can be divided by the monomial $Y_0$, which contradicts the assumption that   $Q$ is of the lowest height  in $\mathcal{S}$.

Now, we finish the proof of Theorem \ref{thm-difference-Holder}.
\end{proof}

\end{document}